\documentclass[11pt]{amsart}
\usepackage{amsmath,amsthm,amsfonts,amssymb,latexsym,epsfig}

\newtheorem{theorem}{Theorem}[section]

\newtheorem{lemma}{Lemma}[section]

\numberwithin{equation}{section}

\theoremstyle{definition}

\theoremstyle{remark}

\begin{document}
\title{A note on Hardy $q$-inequalities}
\author{Peng Gao}
\address{Department of Mathematics, School of Mathematics and System Sciences, Beijing University of Aeronautics and Astronautics, P. R. China}
\email{penggao@buaa.edu.cn}
\subjclass[2000]{Primary 26D15} \keywords{Hardy's inequality}
\thanks{The author is supported in part by NSFC grant 11371043 and the Fundamental Research Funds for the Central Universities.} 

\begin{abstract}We give a simple proof of a recently result concerning Hardy $q$-inequalities. 
\end{abstract}

\maketitle
\section{Introduction}
\label{sec1}
   
   For $p \geq 1$, let $l^p$ be the Banach space of all complex sequences ${\bf x}=(x_n)_{n \geq 1}$ with norm
\begin{equation*}
   \|{\bf x}\|_p: =(\sum_{n=1}^{\infty}|x_n|^p)^{1/p} < \infty.
\end{equation*}

     Let $C=(c_{n,k})$ be a matrix acting on the $l^p$ space, the $l^{p}$ operator norm of $C$ is
   defined as
\begin{equation*}
\label{02}
    \| C \|_{p,p}=\sup_{\|{\bf x}\|_p = 1}\Big | \Big |C {\bf x}\Big | \Big |_p.
\end{equation*}

   The celebrated Hardy's inequality (\cite[Theorem 326]{HLP}) asserts that for $p>1$,
\begin{equation*}
\label{eq:1} \sum^{\infty}_{n=1}\Big{|}\frac {1}{n}
\sum^n_{k=1}x_k\Big{|}^p \leq \Big (\frac
{p}{p-1} \Big )^p\sum^\infty_{n=1}|x_n|^p.
\end{equation*}

    Let $0 < q < 1$ be fixed. The definite $q$-integral or the $q$-Jackson integral (see\cite{MOP}) of a function $f : [0, b) \rightarrow R, 0<b< \infty$, is defined as follows:
\begin{align*}
   \int^x_0f(t)dt=(1-q)x\sum^{\infty}_{k=0}q^kf(q^kx), \quad x\in [0,b).
\end{align*}
   In the theory of $q$-analysis the $q$-analogue $[\alpha]_q$ of a number $\alpha \in R$ is defined by
\begin{align*}
  [\alpha]_q =\frac {1-q^{\alpha}}{1-q}.
\end{align*}

    Recently,  Maligranda, Oinarov and Persson \cite{MOP} considered the $q$-analogues of the classical Hardy type
inequalities, one of their main results is stated as follows: 
\begin{theorem}[{\cite[Theorem 2.3]{MOP}}]
Let $0<q<1$, $\alpha < 1-1/p$. If either $1 \leq p < \infty$ and $f \geq 0$ or $p < 0$ and $f > 0$, then the following inequality
\begin{align*}
   \int^1_0x^{p(\alpha-1)}\Big ( \int^x_0 t^{-\alpha}f(t)d_qt \Big)^pd_qx < \frac 1{[1-1/p-\alpha]^p_q}\int^1_0f^p(t)d_qt
\end{align*}
   holds and the constant $[1-1/p-\alpha]_q$ is best possible.
\end{theorem}

   Using the above result, it is derived in \cite[(27)]{MOP} (here we replace $q^{1-1/p-\alpha}$ in \cite[(27)]{MOP} by $q$) that if $0 < q < 1$  and either $p > 1$ or $p < 0$, then for $a_n>0$,
\begin{align}
\label{1.02}
    \sum^{\infty}_{n=1} \Big ( \frac 1{q^{n}} \sum^{\infty}_{k=n}q^{k}a_k\Big )^p  \leq \frac 1{(1-q)^p}\sum^{\infty}_{n=1}a^p_n,
\end{align}
   where the constant $(1-q)^{-p}$ is best possible. 

   In fact, it is easy to see that Theorem \ref{thm1} is equivalent to inequality \eqref{1.02}. It is therefore our goal in this note to give a simple proof of Theorem \ref{thm1} by proving inequality \eqref{1.02}. We only consider the case $p>1$ here and we prove in the next section the following generalization of inequality \eqref{1.02}:

\begin{theorem}
\label{thm1} Let $p > 1$, let $\{ a_n \}^{\infty}_{n=1}$ be a sequence of non-negative numbers with $a_1>0$ and $S=\sum^{\infty}_{n=1}a_n<\infty$. Let $A=(a_{i,j})$ be an upper triangular matrix satisfying $a_{i,j}=a_{j-i+1}, j \geq i$. Then $\| A \|_{p, p}=S$.
\end{theorem}
\section{Proof of Theorem \ref{thm1}}
\label{sec 2} \setcounter{equation}{0}
    We first note the following generalization of the well-known Schur's test:
\begin{lemma}[{\cite[Theorem 1]{B&J}}]
\label{lem2.1} Let $p>1$ be fixed and let $A=(a_{i,j})$ be a matrix with non-negative entries. If there exist positive numbers $U_1, U_2$ and a positive matrix $(b_{i,j})$, such that
\begin{eqnarray*}
  \sum_{j} a_{i,j}b^{1/p}_{i,j}  &\leq & U_1,  \\
 \sum_{i}a_{i,j}b^{-1/q}_{i,j} &\leq & U_2.
\end{eqnarray*} 
   Then
\begin{equation*}
  \|A \|_{p,p} \leq U^{1/q}_1U^{1/p}_2.
\end{equation*}
\end{lemma}
   
    We now apply the above lemma to the matrix $A$ given in Theorem \ref{thm1} by taking $b_{i,j}=1$ for all $i,j$. Noting that in this case $U_1=S, U_2 \leq S$, we obtain immediately $\|A \|_{p,p} \leq S$.
    
    To see that the constant $S$ is best possible, we note that for any given $\epsilon >0$, there exists an integer $N>0$ such that $\sum^{N}_{i=1}a_{i}>S-\epsilon$. Now consider the sequence ${\bf x}$ with $x_i=1, 1 \leq i \leq M$ for some integer $M>1$ and $x_i=0, i >M$. When $M>N$, the first $M-N+1$ rows of $A{\bf x}$ have row sum $\geq \sum^{N}_{i=1}a_{i}>S-\epsilon$. It follows that $\|A{\bf x} \|_p / \|{\bf x} \|_p >S-2\epsilon$ when $M$ is large enough. As $\epsilon$ is arbitrary, this shows that the constant $S$ is best possible and the proof of Theorem \ref{thm1} is completed.
 


\end{document}